\documentclass[12pt,reqno]{amsart}
\usepackage{eurosym}
\usepackage{amsfonts}
\usepackage{graphicx}
\usepackage{amscd}
\usepackage{amsmath}
\usepackage{epsfig}
\usepackage[usenames]{color}
\usepackage{subfigure}
\usepackage[english]{babel}

\providecommand{\U}[1]{\protect\rule{.1in}{.1in}}
\providecommand{\U}[1]{\protect\rule{.1in}{.1in}}
\allowdisplaybreaks
\newtheorem{theorem}{Theorem}
\newtheorem{lemma}{Lemma}
\newtheorem{proposition}{Proposition}
\theoremstyle{definition}

\theoremstyle{plain}
\newtheorem{acknowledgement}{Acknowledgement}

\newtheorem{corollary}{Corollary}

\newtheorem{definition}{Definition}
\newtheorem{example}{Example}

\numberwithin{equation}{section}


\begin{document}
\title{ NUMERICAL SIMULATIONS FOR FRACTIONAL DIFFERENTIAL EQUATIONS OF HIGHER ORDER AND A WRIGHT-TYPE TRANSFORMATION}
\author{ M. Nacianceno\textsuperscript{1}, T. Oraby\textsuperscript{1}, H.
Rodrigo\textsuperscript{1}, Y. Sepulveda\textsuperscript{1}, J. Sifuentes%
\textsuperscript{1}, E. Suazo\textsuperscript{1, 4}, T. Stuck%
\textsuperscript{3}, J. Williams\textsuperscript{2}, }
\maketitle
\pagestyle{plain} 

\begin{abstract}
In this work, a new relationship is established between the solutions of
higher fractional differential equations and a Wright-type transformation.
Solutions could be interpreted as expected values of functions in a random
time process. As applications, we solve the fractional beam equation,
fractional electric circuits with special functions as external sources, and
derive d'Alembert's formula for the fractional wave equation. Due to this
relationship, we present two methods for simulating solutions of fractional
differential equations. The two approaches use the interpretation of the
Caputo derivative of a function as a Wright-type transformation of the higher derivative of the
function. In the first approach, we use the Runge-Kutta method of hybrid
orders 4 and 5 to solve ordinary differential equations combined with the
Monte Carlo integration to conduct the Wright-type transformation. The
second method uses a feedforward neural network to simulate the fractional
differential equation. 
\end{abstract}

\vspace{5pt}

\begin{enumerate}
\item {\tiny School of Statistical and Mathematical Sciences, The University
of Texas, Rio Grande Valley, 1201 W University Dr, Edinburg, TX 78539 }

\item {\tiny Department of Mathematics, The University of Minnesota, 511 Bruininks Hall 222 Pleasant Street S.E.. Minneapolis MN 55455. }

\item {\tiny Department of Physics, Lamar University, Archer Physics A112
Lamar University 4400 MLK Parkway P.O. Box 10046 Beaumont, Texas 77710 }

\item {\tiny Corresponding Author, erwin .suazo@utrgv.edu.}
\end{enumerate}

\medskip 


\section{Introduction}

While initiated in the latter part of the 18th century by Leibniz, Newton,
and l'Hopital, fractional calculus has garnered significant attention in the
past decades. Fractional differential equations (FDEs) involve fractional derivatives of functions. These equations extend the
concept of ordinary (ODEs) and partial differential equations (PDEs) by allowing the use of
fractional orders of differentiation. Just as differential
equations describe how a function's derivative relates to the function
itself, fractional differential equations describe how fractional
derivatives relate to the function. Solving fractional differential
equations can be challenging due to the nonlocal and noninteger nature of
the derivatives involved. Analytical solutions are not always possible and
numerical methods are often employed to approximate solutions.

Biological systems often exhibit behaviors that involve long-range and memory
interactions \cite{Mainardibook2010}. Fractional differential equations can
be employed to model processes such as enzyme reactions, population dynamics,
and the spread of diseases in a more accurate manner. Fractional calculus
can be used to model long-range memory financial processes, such as
volatility in financial markets. Fractional differential equations can
provide more accurate models for predicting financial fluctuations; for
these and more applications see \cite{Jin2021}, \cite{Meerschaert2012}, \cite%
{Mainardibook2010} and references therein.

In this work, we use an alternative approach to solve fractional ordinary and partial
differential equations of higher order as the expected value of a
random time process $\mathcal{T}_{\beta }(t)$. Let $g_{\beta }(x;t)$ be the
probability density function of the random-time process $\mathcal{T}_{\beta
}(t)$ for all $t>0$ ~\cite{MainConsi2020,Gorenflo2010}.  The random time process $\mathcal{T}_{\beta
}(t)$ can be depicted as the inverse of a $\beta $ stable subordinator with
density $g_{\beta }(x,t)$, see \cite{Bingham1971, Meerschaert2012,
Meerschaert2013}. A subordinator is a non-decreasing L\'{e}vy process, i.e.
a non-negative process with independent and stationary increments. It is
used as a random time, or operational time, in defining time-changed
processes. Its density decays algebraically as $1/t^{\beta +1}$ as $%
t\rightarrow \infty $, with $\beta \in (0,1)$. A different type of time change can be made by using the inverse subordinator (i.e., the hitting
time of a subordinator), which sometimes leads to sub-diffusion processes
\cite{Meerschaert2012}. Using the latter, we apply a numerical approach
based on Monte Carlo integration to simulate solutions of fractional
ordinary and fractional partial differential equations. 

In addition, we present an alternative version of Lemma \ref{lemma1}
presented in \cite{Oraby}, please see Section 3. Lemma \ref{lemma2} tries to establish a connection
between fractional derivatives of higher order and expected values of
derivatives of the same function, see (\ref{integralversion}). However,
condition (\ref{transformation21}) is trivial for $n=1$, but not
trivial for $n>1$. Therefore, Lemma \ref{lemma2} in Section 3 provides a
weaker condition. It was inspired by the work of Dr. Mark Meerschaert and
his work on D'Alembert's formula. Sadly, Dr. Mark Meerschaert passed away
recently, but he contributed tremendously to the theory of fractional
calculus.

This paper is organized as follows: In Section 2, we review important
results of fractional calculus and in particular properties for $g_{\beta
}(x;t)$ and a Wright-type transformation fundamental for our work. In
Section 3, in the Lemma \ref{lemma2} a new relationship is established
between the solutions of higher fractional differential equations and the
Wright-type transformation. Lemma \ref{lemma2} allows us to solve fractional
differential equations of higher orders with certain initial conditions, see
Theorem \ref{Theorem1}. In Theorem \ref{Theorem1} solutions could be
interpreted as expected values of functions in a random-time process. As a
first application, we solve the fractional beam equation. In Section 4, we
continue the applications of Theorem 1, solving fractional electric
circuits with special functions as external sources such as the Mittag-Leffler
function and other special functions. Similar work is also presented for the
homogeneous case. In Section 5, we use the expected value interpretation
presented in Theorem \ref{Theorem1} of the solutions of fractional equations
to perform Monte Carlo simulations of their solution. The idea is to
generate a new but random timeline $L(t)$ using $g_{\beta }(\cdot ,t)$ for
each $t$ and then use any ODE numerical solver, such as Runge-Kutta, to
simulate it on that new timeline and finally find the average of the
solution, see \cite{Oraby}.  In Section 6, we provide a new derivation of d'Alembert's
formula for fractional wave equations as an application of Lemma 2.
 Finally, in Section 7, we show how to use
feedforward neural networks to simulate solutions of
fractional differential equations. The main difference between fractional
and ordinary measurements is in the number of past measurements (memory
length) required to predict the measurement of the next step. We found that
ODEs require one past step, while FDEs require more than one (longer
memory).
\section{On a Wright-type transformation}

In this Section, we review important definitions and classical results that
we will need. Let $D^{n}$ be the Leibniz integer-order differential operator
given by 
\begin{equation*}
D^{n}f=\dfrac{d^{n}f}{dt^{n}}=f^{(n)},
\end{equation*}%
and let $J^{n}$ be an integration operator of integer order given by 
\begin{equation}
J^{n}f(t)=\dfrac{1}{n-1!}\int_{0}^{t}(t-\tau )^{n-1}f(\tau )d\tau ,
\label{eq:intop}
\end{equation}%
where $n\in \mathbb{Z}^{+}$. Let us use $D=D^{1}$ for the first derivative.
For fraction-order integrals, we use 
\begin{equation}
J^{n-\beta }f(t)=\dfrac{1}{\Gamma (n-\beta )}\int_{0}^{t}(t-\tau )^{n-\beta
-1}f(\tau )d\tau ,
\end{equation}%
where $n-1<\beta \leq n$. Now, define the Caputo fractional differential
operator $D_{C}^{\beta }$ to be 
\begin{equation*}
D_{c}^{\beta }f(t)=J^{n-\beta }D^{n}f(t),
\end{equation*}%
where $n-1<\beta \leq n$, for $n\in \mathbb{N}$.

The Wright function is a special function of importance to fractional
calculus and is defined by \cite{MainGloVivo2001}: 
\begin{equation}
W_{\beta ,\alpha }(z)=\sum\limits_{k=0}^{\infty }\dfrac{z^{k}}{k!\Gamma
(\beta k+\alpha )} .  \label{eq:wf}
\end{equation}

$\beta >-1,\alpha \in \mathbb{C},z\in \mathbb{C}$. Also, we will need the
definition of the Mitagg-Leffler function, which is a generalization of the
exponential function.

\begin{definition}
The Mitagg-Leffler function is defined by \textbf{\ } 
\begin{equation*}
E_{\beta ,\alpha }(z)=\sum_{k=0}^{\infty }\frac{z^{k}}{\Gamma (\beta {k}%
+\alpha )},\hspace{5pt}z\in \mathbb{C},\hspace{5pt}\alpha ,\beta >0.
\end{equation*}
\end{definition}

Next, we introduce the main character of this work, the $g_{\beta }(x;t)$,
which is a Wright-type function.

\begin{definition}
$g_{\beta }$ is defined as a Wright-type function 
\begin{equation*}
g_{\beta }(x;t)=\frac{1}{t^{\beta }}W_{-\beta ,1-\beta }\left( \frac{-x}{%
t^{\beta }}\right) =\frac{1}{t^{\beta }}\sum_{k=0}^{\infty }\frac{\left(
-x\right) ^{k}}{t^{\beta k}k!\Gamma (-\beta (k+1)+1)}
\end{equation*}
where $0<\beta \leq 1,\hspace{5pt}t>0,\hspace{5pt}s\geq 0.$ $g_{\beta }(x;t)$
is a probability density function of the random time process $\mathcal{T}%
_{\beta }(t)$ for all $t>0$ ~\cite{MainConsi2020,Gorenflo2010}. $\mathcal{T}%
_{\beta }(t)$ can be seen as the inverse of a $\beta $ stable subordinator
with density $g_{\beta }(x,t)$, see \cite{Bingham1971, Meerschaert2012,
Meerschaert2013}.
\end{definition}

We use $g_{\beta }(x;t)$ to define the following transform.

\begin{definition}
We define a Wright-type transform for a function $f(t)$ as 
\begin{equation*}
f_{\beta }(t)=\int_{0}^{\infty }g_{\beta }(s;t)f(s)ds
\end{equation*}%
where $t>0,\hspace{5pt}0<\beta \leq 1.$ We will use the alternative notation 
$f_{\beta }(t)=$ $(f(t))_{\beta }.$
\end{definition}

The following results are properties of $g_{\beta }(x;t)$. Some of them are
well known, see \cite{Mainardibook2010}, and others have been recently
proved in \cite{Oraby}. We set these properties as a proposition for the
convenience of the reader, and they will play a fundamental role in our
numerical simulations.

\begin{proposition}
\label{proposition1} The following properties for $g_{\beta }$, $0<\beta<1$, are valid:
\end{proposition}

\begin{enumerate}
\item $g_{\beta }$ is a probability distribution%
\begin{equation*}
\int_{0}^{\infty }g_{\beta }(s;t)ds=1.
\end{equation*}

\item $g_{\beta }$ has the following Laplace transforms 
\begin{equation}
\mathcal{L}\left( g_{\beta }(\cdot ;t)\right) (s)=\int_{0}^{\infty
}e^{-sx}g_{\beta }(x;t)dx=E_{\beta }\left( -st^{\beta }\right)  \label{eq:lw}
\end{equation}%
for $\Re (s)>0$ and moments $\mathbb{E}\left[ (\mathcal{T}_{\beta }(t))^{k}%
\right] =\Gamma (k+1)\dfrac{t^{k\beta }}{\Gamma (k\beta +1)}$ for $k\geq 1$~%
\cite{Piryatinska2005,Kumar2018}. At the same time 
\begin{equation}
\int_{0}^{\infty }e^{-st}g_{\beta }(x;t)dt=s^{\beta -1}e^{-xs^{\beta }}.
\label{eq:lw2}
\end{equation}%

\item $g_{\beta }$ has the following moments:
\begin{equation*}
\int_{0}^{\infty }s^{n}g_{\beta }(s;t)ds=E[S^{n}]=n!\frac{t^{n\beta }}{%
\Gamma (n\beta +1)}.
\end{equation*}

\item The following \textbf{Taylor Expansion for }$f_{\beta }(t)$ holds:
\textbf{\ } 
\begin{equation*}
f_{\beta }(t)=\int_{0}^{\infty }f(s)g_{\beta }(s;t)ds=\int_{0}^{\infty
}\sum_{n=0}^{\infty }f^{(n)}(0)\frac{s^{n}}{n!}g_{\beta
}(s;t)ds=\sum_{n=0}^{\infty }f^{(n)}(0)\frac{t^{n\beta }}{\Gamma (n\beta +1)}%
.
\end{equation*}
\end{enumerate}

From the proposition above we can make a table of the Wright beta transform
for different functions. This table will be useful in Section 4 where we will
solve nonhomogeneous fractional ODEs.
\begin{table}
\begin{tabular}{|l|l|}
    \hline
    $f(t)$ & $f_{\beta }(t)$ \\ \hline
    $t$ & $\frac{t^{\beta }}{\Gamma (\beta +1)}$ \\ \hline
    $t^{2}$ & $\frac{2t^{2\beta }}{\Gamma (2\beta +1)}$ \\ \hline
    $t^{n}$ & $\frac{n!t^{n\beta }}{\Gamma (n\beta +1)}$ \\ \hline
    $\cos (t)$ &  $\cos_{\beta}(t)=\sum_{n=0}^{\infty }\frac{(-1)^{n} t^{2n\beta}}{\Gamma
    (2n\beta +1)}=E_{2\beta}(-t^{2\beta})$ \\ \hline
    $\sin (t)$ & $\sin_{\beta}(t)=\sum_{n=0}^{\infty }\frac{(-1)^{n}t^{\beta (2n+1)}}{\Gamma
    ((2n+1)\beta +1)}=t^{\beta}E_{2\beta,\beta+1}(-t^{2\beta})$ \\ \hline
    $e^{t}$ & $\sum_{n=0}^{\infty}\frac{t^{n\beta }}{\Gamma (n\beta +1)}=E_{\beta}(t^{\beta})$ \\ \hline
    $\cosh (t)$ &  $\cosh_{\beta}(t)=\sum_{n=0}^{\infty} \frac{t^{(2n)\beta}}{\Gamma(2n\beta+1)}=E_{2\beta}(t^{2\beta})$ \\ 
    \hline
    $\sinh (t)$ &  $\sinh_{\beta}(t)=\sum_{n=0}^{\infty} \frac{t^{(2n+1)\beta}}{\Gamma((2n+1)\beta+1)}=t^{\beta}E_{2\beta,\beta+1}(t^{2\beta})$ \\ \hline
    $e^{t}\cos (t)$ & $\sum_{n=0}^{\infty} [\sum_{k=0}^{n} \binom{n}{2k} {(-1)^{k}}] \frac{t^{\beta n}}{\Gamma
    (n\beta +1)}$ \\ \hline
    $e^{-t}\cos (t)$ & $\sum_{n=0}^{\infty} [\sum_{k=0}^{n} \binom{n}{2k} {(-1)^{n-k}}]\frac{t^{\beta n}}{\Gamma
    (n\beta +1)}$  \\ \hline
    $e^{t}\sin (t)$ & $\sum_{n=0}^{\infty} [\sum_{k=0}^{n} \binom{n}{2k+1} {(-1)^{k}}]\frac{t^{\beta n}}{\Gamma
    (n\beta +1)}$  \\ \hline
    $e^{-t}\sin (t)$ & $\sum_{n=0}^{\infty} [\sum_{k=0}^{n} \binom{n}{2k+1} {(-1)^{n-k-1}}]\frac{t^{\beta n}}{\Gamma
    (n\beta +1)}$  \\ \hline
    $g(t)\cos (t)$ & $\sum_{n=0}^{\infty} [\sum_{k=0}^{n} \binom{n}{2k} {(-1)^{k}} {g^{(n-2k)}(0)}]\frac{t^{\beta n}}{\Gamma
    (n\beta +1)}$\\ \hline
    $g(t)\sin (t)$ & $\sum_{n=0}^{\infty} [\sum_{k=0}^{n} \binom{n}{2k+1} {(-1)^{k}} {g^{(n-2k-1)}(0)}]\frac{t^{\beta n}}{\Gamma
    (n\beta +1)}$ \\ \hline
    $g(t)\cosh (t)$ & $\sum_{n=0}^{\infty} [\sum_{k=0}^{n} \binom{n}{2k} {g^{(n-2k)}(0)}]\frac{t^{\beta n}}{\Gamma
    (n\beta +1)}$ \\ \hline
    $g(t)\sinh (t)$ & $\sum_{n=0}^{\infty} [\sum_{k=0}^{n} \binom{n}{2k+1} {g^{(n-2k-1)}(0)}]\frac{t^{\beta n}}{\Gamma
    (n\beta +1)}$ \\ \hline 
    \end{tabular}
    \caption{ A Wright beta transform
for different functions. This table will be useful in Section 4 where we will
solve nonhomogeneous fractional ODEs}
    \end{table}

\section{Higher order fractional ordinary differential equations}

In this Section, we consider a new relationship between solutions of higher
fractional differential equations and the Wright-type transformation, thus improving the following Lemma presented previously in \cite{Oraby}.
\begin{lemma}
\label{lemma1} Let $f\in C^{n}([0,\infty ))$ such that 
\begin{equation}
\mathcal{L}\left( f_{\beta }^{(n-1)}(\cdot )\right) (s)=s^{\beta -1}\mathcal{%
L}\left( f^{(n-1)}(\cdot )\right) \left( s^{\beta }\right)
\label{transformation21}
\end{equation}%
for $s\in \lbrack 0,\infty )$ and $f_{\beta }^{(n-1)}(0)=f^{(n-1)}(0)$. For $%
0<\beta \leq 1$, the following holds 
\begin{equation}
D_{c}^{n+\beta -1}f_{\beta }(t)=\int_{0}^{\infty }f^{(n)}(x)g_{\beta }(x;t)dx
\label{integralversion}
\end{equation}%
for $n\geq 1$.
\end{lemma}

We point out the following lemma, which characterizes higher fractional
derivatives of different orders in terms of Wright-type transformations; see
(\ref{ncase}), which improves the presentation given in Lemma \ref{lemma1}, see also \cite{Oraby}. In
particular, this Lemma does not require the restricted condition (\ref%
{transformation21}).

\begin{lemma}
\label{lemma2} Under the assumption of $f_{\beta }(0)=f(0),$ $D_{c}^{\beta
}f_{\beta }(0)=f^{\prime }(0),$ $D_{c}^{\beta }D_{c}^{\beta }f_{\beta
}(0)=f^{^{\prime \prime }}(0)...$ and $\underbrace{D_{c}^{\beta }\cdot \cdot
\cdot D_{c}^{\beta }D_{c}^{\beta }}_{\text{$n$ times}}f_{\beta
}(0)=f^{(n)}(0),$ the following holds:
\begin{equation}
\underbrace{D_{c}^{\beta }\cdot \cdot \cdot D_{c}^{\beta }D_{c}^{\beta }}_{%
\text{$n$ times}}f_{\beta }=\int_{0}^{\infty }f^{(n)}(x)g_{\beta }(x,t)dx.
\label{ncase}
\end{equation}
\end{lemma}

\textbf{Proof}

We will prove that both parts of the equation have the same Laplace
transform. The case $n=1$ has been proved previously  in \cite{Oraby}$.$
Indeed $D_{c}^{\beta }f_{\beta }(t)=\int_{0}^{\infty }f^{\prime }(x)g_{\beta
}(x;t)dx,$ where $0<{\beta }<1.$ Taking the Laplace transform on both sides
of (\ref{ncase}), on the left-hand side we have 
\begin{eqnarray*}
\mathcal{L}\{\underbrace{D_{c}^{\beta }\cdot \cdot \cdot D_{c}^{\beta
}D_{c}^{\beta }}_{\text{$n$ times}}f_{\beta }\} &=&s^{\beta }\overline{%
\underbrace{D_{c}^{\beta }\cdot \cdot \cdot D_{c}^{\beta }D_{c}^{\beta }}_{%
\text{$n-1$ times}}f_{\beta }(t)}-s^{\beta -1}\underbrace{D_{c}^{\beta
}\cdot \cdot \cdot D_{c}^{\beta }D_{c}^{\beta }}_{\text{$n-1$ times}%
}f_{\beta }(0) \\
&=&s^{\beta }(s^{n\beta -1}\overline{f}(s^{\beta })-s^{(n-1)\beta
-1}f_{\beta }(0)-s^{(n-2)\beta -1}D_{c}^{\beta }f_{\beta }(0)-...) \\
&&-s^{\beta -1}\underbrace{D_{c}^{\beta }\cdot \cdot \cdot D_{c}^{\beta
}D_{c}^{\beta }}_{\text{$n-1$ times}}f_{\beta }(0) \\
&=&s^{(n+1)\beta -1}\overline{f}(s^{\beta })-s^{(n)\beta -1}f_{\beta
}(0)-s^{n\beta -\beta -1}D_{c}^{\beta }f^{{}}(0)...-s^{\beta -1}\underbrace{%
D_{c}^{\beta }\cdot \cdot \cdot D_{c}^{\beta }D_{c}^{\beta }}_{\text{$n-1$
times}}f_{\beta }(0).
\end{eqnarray*}%
On the other hand, taking the Laplace transform on the right-hand side\ and
using our favorite equation (\ref{eq:lw}), we have 
\begin{eqnarray*}
\mathcal{L}\left\{ {\int_{0}^{\infty }f^{(n)}(x)g_{\beta }(x;t)dx}\right\}
&=&s^{\beta -1}\mathcal{L}\{{f}^{(n)}{(\cdot )}\}(s^{\beta }) \\
&=&s^{\beta -1}\left( s^{n\beta }\overline{f}(s^{\beta })-s^{(n-1)\beta
}f(0)-s^{(n-2)\beta }f^{\prime }(0)...-f^{(n-1)}(0)\right) \\
&=&s^{(n+1)\beta -1}\overline{f}(s^{\beta })-s^{(n)\beta -1}f(0)-s^{n\beta
-\beta -1}f^{\prime }(0)...-s^{\beta -1}f^{(n-1)}(0).
\end{eqnarray*}%
Therefore under the assumptions of $f_{\beta }(0)=f(0),$ $D_{c}^{\beta
}f_{\beta }(0)=f^{\prime }(0),$ $D_{c}^{\beta }D_{c}^{\beta }f_{\beta
}(0)=f^{^{\prime \prime }}(0)...$ and $\underbrace{D_{c}^{\beta }\cdot \cdot
\cdot D_{c}^{\beta }D_{c}^{\beta }}_{\text{$n$ times}}f_{\beta
}(0)=f^{(n)}(0)$, we obtain the equality$.$

Thanks to Lemma 2, we obtain the following Theorem which is a significant
improvement to Lemma 1. In the following Theorem we don't require condition (\ref%
{transformation21}) which seems to be quite restricted and was presented
in \cite{Oraby}.

\begin{theorem}
\label{Theorem1} \textbf{\ }For non-homogeneous fractional differential
equations of the form: 
\begin{equation*}
\sum_{k=1}^{n-1}a_{k}\underbrace{D_{c}^{\beta }\cdot \cdot \cdot
D_{c}^{\beta }D_{c}^{\beta }}_{k\text{ times}}D_{C}^{\beta
}y(t)+a_{n+1}y(t)=F_{\beta }(t)
\end{equation*}%
such that $f_{\beta }(0)=z(0),$ $D_{c}^{\beta }f_{\beta }(0)=z^{\prime }(0),$
$D_{c}^{\beta }D_{c}^{\beta }f_{\beta }(0)=z^{^{\prime \prime }}(0)...$ and $%
\underbrace{D_{c}^{\beta }\cdot \cdot \cdot D_{c}^{\beta }D_{c}^{\beta }}_{%
\text{$n$ times}}f_{\beta }(0)=z^{(n)}(0)$ for $0<\beta \leq 1$, the
solution is given by 
\begin{equation*}
y_{\beta }(t)=\int_{0}^{\infty }z(s)g_{\beta }(s;t)ds=\mathbb{E}\left[ z(%
\mathcal{T}_{\beta }(t))\right],
\end{equation*}%
where the process $\mathcal{T}_{\beta }(t)$ can be seen as the inverse of a $%
\beta $ stable subordinator whose density is $g_{\beta }(x,t)$, $F_{\beta
}(t)=\int_{0}^{\infty }F(s)g_{\beta }(s;t)ds$, and $z(t)$ is the solution of
the linear ordinary differential equation with the same initial conditions. 
\begin{equation*}
\sum_{k=1}^{n}a_{k}z^{(k)}(t)+a_{n+1}z(t)=F(t)
\end{equation*}

$z^{n-1}(0)=z_{n},...z(0)=z_{0}.$
\end{theorem}

The following corollary solves second-order homogeneous fractional ODEs.

\begin{corollary}
\label{Corollary1} Consider the equation 
\begin{equation*}
a_{2}D_{c}^{\beta }D_{c}^{\beta }y_{\beta }(t)+a_{1}D_{c}^{\beta }y_{\beta
}(t)+a_{0}y_{\beta }(t)=0,
\end{equation*}
if $r_1$ and $r_2$ are the roots of the characteristic polynomial associated to the analogous second order differential equation, then the solution has the form:\\
Case 1: Distinct and Real Roots 
\begin{equation*}
\sum_{n=0}^{\infty }{[c_{2}r_{2}^{n}+c_{1}r_{1}^{n}]\frac{t^{n\beta }}{%
\Gamma {(n\beta +1)}}}=c_{2}E_{\beta }{(r_{2}t^{\beta })}+c_{1}E_{\beta }{%
(r_{1}t^{\beta }).}
\end{equation*}

Case 2: Repeating Roots 
\begin{equation*}
\sum_{n=0}^{\infty }{[c_{2}r^{n}+c_{1}nr^{n-1}]\frac{t^{n\beta }}{\Gamma {%
(n\beta +1)}}}=c_{2}E_{\beta }(rt^{\beta })+c_{1}\frac{\partial }{\partial r}%
E_{\beta }(rt^{\beta }).
\end{equation*}

Case 3: Complex Roots Let $r_{2}=\lambda +\mu i,r_{1}=\lambda -\mu i $ 
\begin{equation*}
\sum_{n=0}^{\infty }{[c_{2}r_{2}^{n}+c_{1}r_{1}^{n}}]\frac{t^{n\beta }}{%
\Gamma {(n\beta +1)}}=c_{2}E_{\beta }{(r_{2}t^{\beta })}+c_{1}E_{\beta }{%
(r_{1}t^{\beta }).}
\end{equation*}
\end{corollary}

\subsection{Application to Fractional Embedded Beam}
Euler studied as an eigenvalue problem how a thin elastic column buckles under a compressive axial. The differential equation
governing the deflection $\phi (t)$ of a thin elastic column subject to a
constant compressive axial force (or load) $P$ when $P$ is applied to its
top is given by%
\begin{equation}
\frac{d^{2}}{dt^{2}}\left( EI\frac{d^{2}\phi (t)}{dt^{2}}\right) +P\frac{%
d^{2}\phi (t)}{dt^{2}}=0  \label{Beam}
\end{equation}

where $E$ is Young's modulus of elasticity and $I$ is the moment of inertia
of a cross-section about a vertical line through a centroid, see \cite{D. Zill}. If we assume
that the column is uniform, $EI$ is a constant. In the next example, we
solve the equivalent fractional high-order differential equation.

Consider a beam of length $L$. As a first application, we solve the fractional embedded beam equation%
\begin{equation*}
EI\,D_{c}^{\beta }D_{c}^{\beta }D_{c}^{\beta }D_{c}^{\beta }\phi (t)=w_{0},
\end{equation*}%
subject to $\phi (0)=0,\phi (L)=0,D_{c}^{\beta }\phi (0)=0,D_{c}^{\beta
}\phi (L)=0.$ Using Theorem \ref{Theorem1} we obtain%
\begin{eqnarray*}
\phi _{\beta }(t) &=&\frac{w_{0}}{24EI}\int_{0}^{\infty
}(x^{4}-2x^{3}L+L^{2}x^{2})g_{\beta }(x;t)dx \\
&=&\frac{w_{0}}{24EI}\left( 4!\frac{t^{4\beta }}{\Gamma (4\beta +1)}-12\frac{%
t^{3\beta }L}{\Gamma (3\beta +1)}+L^{2}\frac{t^{2\beta }L}{\Gamma (2\beta +1)%
}\right) .
\end{eqnarray*}

\begin{figure}[h]
{\includegraphics{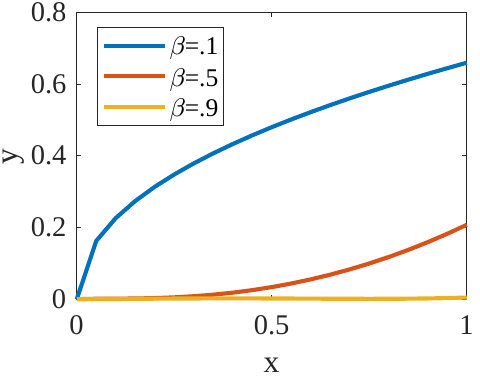}}
\caption{Monte Carlo simulation of the solution of $EI\,D_{c}^{\beta }D_{c}^{\beta }D_{c}^{\beta }D_{c}^{\beta }\phi (t)=w_{0}$ with $L=1$ and $%
EI=w_0 $.}
\label{fig1_1}
\end{figure}
Figure 1 shows a Monte Carlo simulation of the solution, for more details see Section 5.
As a second application we solve, the following equivalent fractional equation corresponding to (\ref%
{Beam})%
\begin{equation*}
EID_{c}^{\beta }D_{c}^{\beta }D_{c}^{\beta }D_{c}^{\beta }\phi
(t)+PD_{c}^{\beta }D_{c}^{\beta }\phi (t)=0,
\end{equation*}%
subject to $\phi (0)=0,\phi (L)=0,D_{c}^{\beta }D_{c}^{\beta }\phi (0)=0$
and $D_{c}^{\beta }D_{c}^{\beta }\phi (L)=0.$ By Theorem 1, the solution is
given by 
\begin{eqnarray*}
\phi _{\beta }(t) &=&\int_{0}^{\infty }(c_{1}\cos (n\pi x/L)+c_{2}\sin
(\alpha x)+c_{3}x+c_{4})g_{\beta }(x;t)dx \\
&=&c_{1}cos_{\beta }((n\pi /L)^{\frac{1}{\beta }}t)+c_{2}sin_{\beta }((n\pi
/L)^{\frac{1}{\beta }}t)+\frac{t^{\beta }L}{\Gamma (\beta +1)}+c_{4} \\
&=&c_{1}E_{2\beta }(-((n\pi /L)^{\frac{1}{\beta }}t)^{2\beta })+c_{2}(n\pi
/L)t^{\beta }E_{2\beta ,\beta +1}(-(n\pi /L)^{2}t^{2\beta })+\frac{t^{\beta
}L}{\Gamma (\beta +1)}+c_{4}.
\end{eqnarray*}

\section{Solving Fractional Electric Circuits}

In this Section, we apply Theorem 1 to solve extensions of fractional
electric circuits with an external source, including an external source
coming from the Mittagg-Leffler function and other special functions;
similar work was presented for the homogeneous case, \cite{Alsaedi}. The
basic equations of electric circuits involving resistors, capacitors, and
inductors are \cite{D. Zill}:
\begin{eqnarray}
I^{\prime \prime }(t)+\frac{1}{LC}I(t) &=&0,  \label{LC} \\
V^{\prime }(t)+\frac{1}{RC}V(t) &=&0,  \label{RC} \\
I^{\prime }(t)+\frac{R}{L}I(t) &=&\frac{V}{L},  \label{RL}
\end{eqnarray}

where the $I(t)$ is the current in the circuit at time $t,$ $L$ is the
inductance, $C$ is the capacitance, $R$ is the resistance and $V$ is the
voltage drop across the circuit. Equation (\ref{LC}) represents the LC
(inductor -capacitor) circuit, equation (\ref{RC}) represents RC (resistor-
capacitor) circuit.

\begin{example}[RC circuit]
\label{exm:exm2} By Theorem 1, the solution of the fractional circuit $%
D_{c}^{\beta }V(t)+V(t)/RC=0$ with $V(0)=V_{0}$ is given by 
\begin{equation*}
V_{\beta }(t)=\int_{0}^{\infty }V_{0}e^{-x/RC}g_{\beta
}(x;t)dx=V_{0}E_{\beta }\left( -t^{\beta }/RC\right).
\end{equation*}%
This solution agrees with the literature \cite{Jin2021}. Furthermore by Theorem 1, $V_{\beta
}(t)=\mathbb{E}[y_{0}e^{\lambda \mathcal{T}_{\beta }(t)}]$.
\end{example}

\begin{example}[LC Circuits]
\label{exm:exm3} The solution of the fractional LC circuit $D_{c}^{\beta
}D_{c}^{\beta }V(t)+V(t)/LC=0$ with $V(0)=0$ and $V^{\prime }(0)=1$, where $%
\omega =1/LC$ is a real-valued constant is given by%
\begin{eqnarray*}
V_{s}(t) &=&\int_{0}^{\infty }g_{\beta }(x,t)\sin (\omega
x)\,dx=\sum_{n=0}^{\infty }\frac{(-1)^{n}(\omega )^{2n+1}}{(2n+1)!}%
\int_{0}^{\infty }g_{\beta }(x,t)\frac{x^{2n+1}}{(2n+1)!}\,dx \\
&=&\sum_{n=0}^{\infty }\frac{(-1)^{n}(\omega ^{\frac{1}{\beta }}t)^{\beta
(2n+1)}}{\Gamma ((2n+1)\beta +1)}=\omega t^{\beta }E_{2\beta ,\beta
+1}(-\omega ^{2}t^{2\beta })=sin_{\beta }(\omega ^{\frac{1}{\beta }}t).
\end{eqnarray*}%
We have applied Theorem 1 to $z(x)=\sin (\omega x)$ which solves $%
D^{2}z(x)=-\omega ^{2}z(x)$, with $z(0)=0$. That is, $y(t)=\mathbb{E}[\sin
(\omega \mathcal{T}_{\beta }(t))]$. This agrees with A.A. Stanislasky's
finding, see equation \cite{A. A. Stanislavsky}. Furthermore, the solution of
the frational LC circuit $D_{C}^{\beta }D_{C}^{\beta }V(t)+V(t)/LC=0$ with $%
y(0)=1$ and $y^{\prime }(0)=0$ where $\omega =1/LC$ is a real-valued
constant is given by 
\begin{eqnarray*}
V_{c}(t) &=&\int_{0}^{\infty }g_{\beta }(x,t)\cos (\omega
x)\,dx=\sum_{k=0}^{\infty }\frac{(-1)^{k}(\omega )^{2k}}{(2k)!}%
\int_{0}^{\infty }g_{\beta }(x,t)\frac{x^{2k}}{(2k)!}\,dx \\
&=&\sum_{k=0}^{\infty }\frac{(-1)^{n}\omega ^{2k}}{(2k)!}\frac{t^{\beta
(2k)}(2k)!}{\Gamma (2k\beta +1)}=\sum_{k=0}^{\infty }\frac{(-1)^{k}(\omega ^{%
\frac{1}{\beta }}t)^{\beta 2k}}{\Gamma (2k\beta +1)} \\
&=&E_{2\beta }(-(\omega ^{\frac{1}{\beta }}t)^{2\beta })=cos_{\beta }(\omega
^{\frac{1}{\beta }}t).
\end{eqnarray*}
\end{example}

This also agrees with A.A. Stanislasky's finding $\mathbf{E}_{2\beta
,1}(\omega ^{2}t^{2\beta })$ \cite{A. A. Stanislavsky}$.$ 

Next, we consider examples of fractional electric circuits with an external
source.

\begin{example}
Solve $D_{c}^{\beta }y_{\beta }(t)+ay_{\beta }(t)=F_{\beta }(t)=E_{\beta
}(t^{\beta }),$ $0<\beta <1$ for the case that $F(t)=e^{t}$. Combining the
solution for the regular ODE $y^{\prime }+ay=F(t)$ and by Theorem 1, we have%
\begin{equation*}
y_{\beta }(t)=\int_{0}^{\infty }\left( \frac{1}{a+1}(e^{s}-e^{-as})+e^{-as}%
\right) g_{\beta }(s;t)ds,
\end{equation*}
\end{example}

and by propostion 1, we obtain 
\begin{equation*}
y_{\beta }(t)=\frac{1}{a+1}E_{\beta }(t^{\beta })+\frac{1}{a+1}E_{\beta
}(-at^{\beta })+y_{0}E_{\beta }(-at^{\beta }).
\end{equation*}
\begin{figure}[h]
{\includegraphics{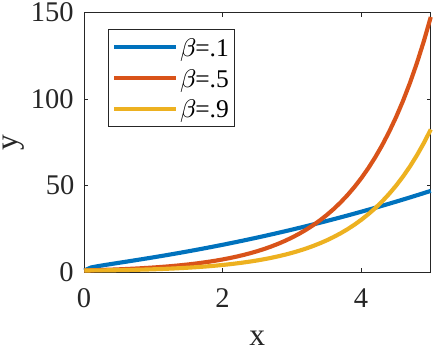}}
\caption{A Monte Carlo simulation of the solution of $D_{c}^{\beta }y_{\beta }(t)+ay_{\beta }(t)=E_{\beta
}(t^{\beta }),$ with $a=1$.}
\label{fig4_1}
\end{figure}
Figure 2 shows a Monte Carlo simulation of the solution; for more details see Section 5.
\begin{example}
Let's consider the following non-homogeneous fractional equation%
\begin{equation*}
D_{c}^{\beta }D_{c}^{\beta }y_{\beta }(t)+\omega ^{2}y_{\beta }(t)=\left(
t^{2}\right) _{\beta }=\frac{2t^{2\beta }}{\Gamma (2\beta +1)},0<\beta <1.
\end{equation*}
\end{example}

Using Theorem 1, we apply the Wright-type transform to the solution of the
standard ODE $y^{^{\prime \prime }}(t)+\omega ^{2}y(t)=t^{2},$ obtaining 
\begin{align*}
y_{\beta }(t)& =\int_{0}^{\infty }\left[ c_{1}\cos (\omega s)+c_{2}\sin
(\omega s)+\left( \frac{s^{2}}{\omega ^{2}}-\frac{2}{\omega ^{4}}\right) %
\right] g_{\beta }(s;t)ds \\
& =c_{1}E_{2\beta }(-(\omega t)^{2\beta })+c_{2}(\omega t)^{\beta }E_{2\beta
,\beta +1}(-(\omega t)^{2\beta })+\frac{1}{\omega ^{2}}\left[ \frac{%
2t^{2\beta }}{\Gamma (2\beta +1)}\right] -\frac{2}{\omega ^{4}}.
\end{align*}

\begin{figure}[h]
{\includegraphics{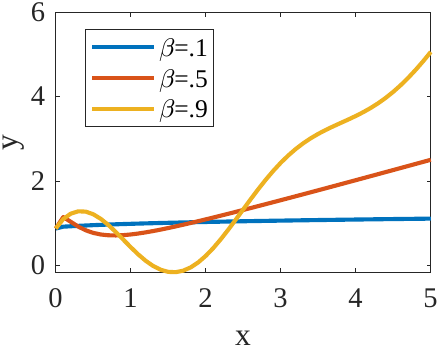}}
\caption{A Monte Carlo simulation of the solution of 
$D_{c}^{\beta }D_{c}^{\beta }y_{\beta }(t)+\omega ^{2}y_{\beta }(t)=\frac{2t^{2\beta }}{\Gamma (2\beta +1)},0<\beta <1$, and $c_1=c_2=1$. }
\label{fig3_1}
\end{figure}
Figure 3 shows a Monte Carlo simulation of the solution; for more details see Section 5.
\begin{example}
$D_{c}^{\beta }D_{c}^{\beta }$\textbf{$y(t)+\omega ^{2}y(t)=t^{\beta
}E_{2\beta ,\beta +1}(-t^{2\beta })=sin_{\beta }(t),$ $0<\beta <1$}
\end{example}

Using Theorem 1, we apply the Wright-type transform to the solution of the
standard ODE $y^{^{\prime \prime }}(t)+\omega ^{2}y(t)=\sin {(t)},$
obtaining 
\begin{align*}
y_{\beta }(t)& =\int_{0}^{\infty }C_{1}\cos {(\omega x)}+C_{2}\sin {(\omega
x)}+\frac{\sin {(x)}}{\omega ^{2}-1}g_{\beta }(x;t)dx \\
& =\int_{0}^{\infty }C_{1}\cos {(\omega x)}g_{\beta
}(x;t)dx+\int_{0}^{\infty }C_{2}\sin {(\omega x)}g_{\beta
}(x;t)dx+\int_{0}^{\infty }\frac{\sin {(x)}}{\omega ^{2}-1}g_{\beta }(x;t)dx
\\
& =C_{1}\cos _{\beta }(\omega t)+C_{2}\sin _{\beta }(\omega t)+\frac{\sin
_{\beta }(t)}{\omega ^{2}-1} \\
& =C_{1}E_{2\beta }(-\omega ^{2}t^{2\beta })+C_{2}\omega t^{\beta }E_{2\beta
,\beta +1}(-\omega ^{2}t^{2\beta })+\frac{t^{\beta }E_{2\beta ,\beta
+1}(-t^{2\beta })}{\omega ^{2}-1}.
\end{align*}

\begin{figure}[h]
{\includegraphics{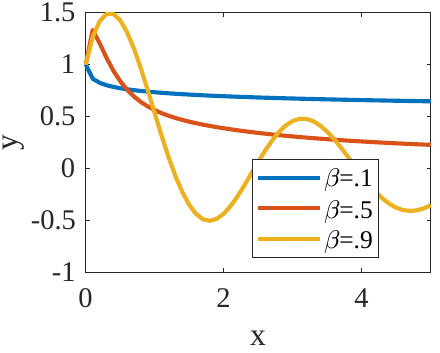}}
\caption{A Monte Carlo simulation of the solution of $D_{c}^{\beta }D_{c}^{\beta }$\textbf{$y(t)+\omega ^{2}y(t)=t^{\beta
}E_{2\beta ,\beta +1}(-t^{2\beta })=sin_{\beta }(t),$ $0<\beta <1$} with $\protect\omega%
=2 $, and $c_1=c_2=1$. }
\label{fig1_1}
\end{figure}
Figure 4 shows a Monte Carlo simulation of the solution; for more details see Section 5
\section{Monte Carlo Simulations of Fractional Differential Equations}

In this section, we will show how to simulate solutions of fractional
differential equations using Monte Carlo methods. The following algorithm
was introduced for first derivatives in \cite{Oraby}. In this work, we
extend these ideas to higher derivatives. The idea is to generate a new but
random timeline $L(t)$ using $g_{\beta }(\cdot ,t)$ for each $t$ and then
use any ODE numerical solver, like Runge-Kutta, to simulate it over that new
timeline and finally find the average of the solution, see \cite{Oraby}.

\textbf{Monte-Carlo Algorithm}

\begin{enumerate}
\item \textbf{Step 1:} For each $t_i$, where $i=1,2,\ldots,n$, generate
random numbers $s_1,\ldots,s_m$ from $g_{\beta}(\cdot,t_i)$.

\item \textbf{Step 2:} Discretize the timeline $L(t)$ from $%
\min(s_1,\ldots,s_m)$ to $\max(s_1,\ldots,s_m)$.

\item \textbf{Step 3:} Solve the ODE using Runge-Kutta from $%
\min(s_1,\ldots,s_m)$ to $\max(s_1,\ldots,s_m)$ and interpolate the solution
at the rest of $s_1,\ldots,s_m$.

\item \textbf{Step 4:} Find the average of the solution over those $%
s_1,\ldots,s_m$.
\end{enumerate}

In the following three examples, we simulate some fractional differential
equations using the Monte Carlo algorithm.

\begin{example}
Let $0<\beta <1$, and consider 
\begin{equation*}
D_{c}^{\beta }y(t)+y(t)=\sin _{\beta }{(t)}
\end{equation*}%
such that $y(0)=0$. It is simulated by solving $z^{\prime }(t)+z(t)=\sin {(t)%
}$ such that $z(0)=0$ using Runge-Kutta of 4 and 5 over simulated timelines $%
L(t)$ for each $t$. See Figure 5.

\begin{figure}[h]
{%
\includegraphics[scale=0.29]{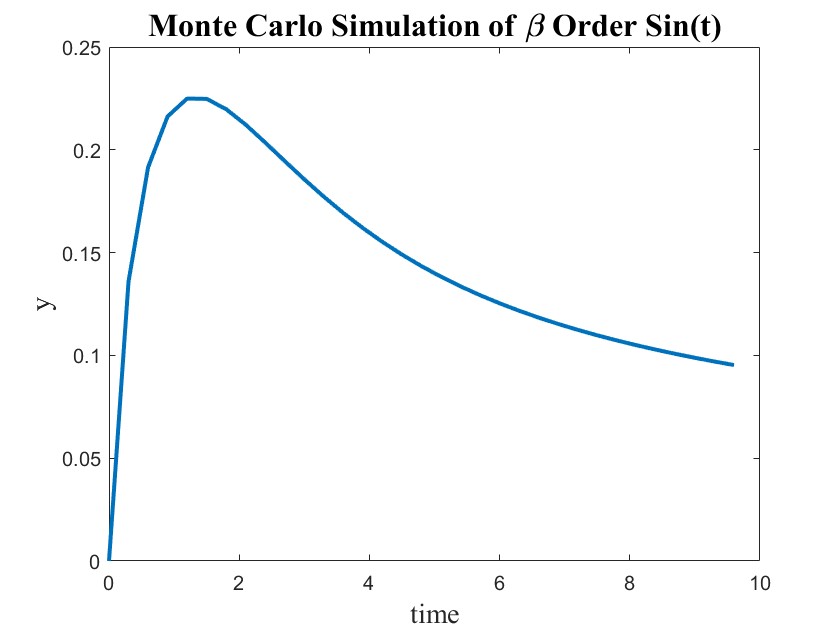}%
} {%
\includegraphics[scale=0.29]{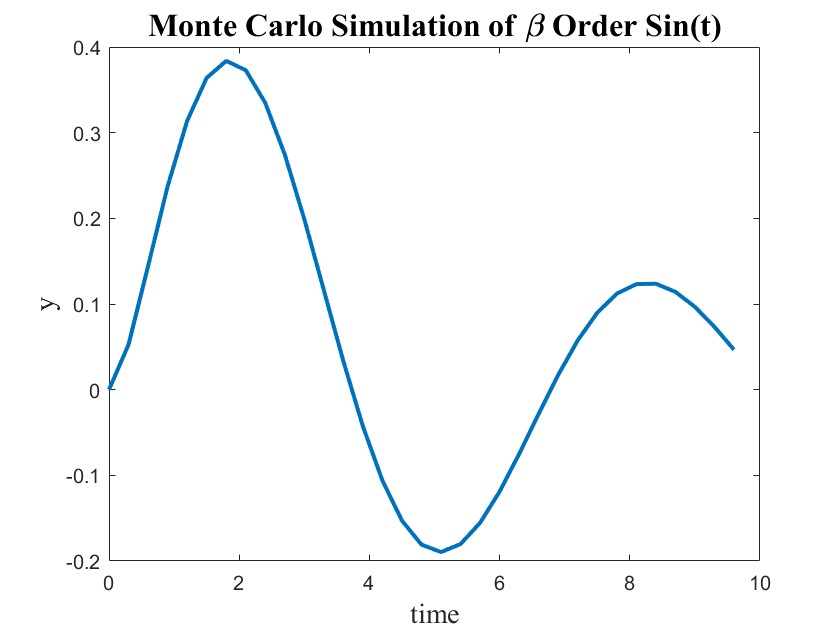}%
}
\caption{A Monte Carlo simulation of the solution of $D_{C}^{\protect\beta %
}y(t)+y(t)=\sin _{\protect\beta }{(t)}$ with $y(0)=0$ evaluated at $\protect%
\beta =0.5$ (left) and $\protect\beta =0.9$ (right).}
\label{fig3}
\end{figure}
\end{example}

\begin{example}
Let $0<\beta <1$, and consider%
\begin{equation*}
D_{c}^{\beta }D_{c}^{\beta }y(t)+\omega ^{2}y(t)=\sin _{\beta }{(t)}
\end{equation*}%
such that $y(0)=0$ and $y^{\prime }(0)=0$. It is simulated by solving $%
z^{\prime \prime }(t)+w^{2}z(t)=\sin {(t)}$ such that $z(0)=0$ and $%
z^{\prime }(0)=0$ using Runge-Kutta of 4 and 5 over simulated timelines $%
L(t) $ for each $t$. See Figure \ref{fig4}.

\begin{figure}[h]
\includegraphics{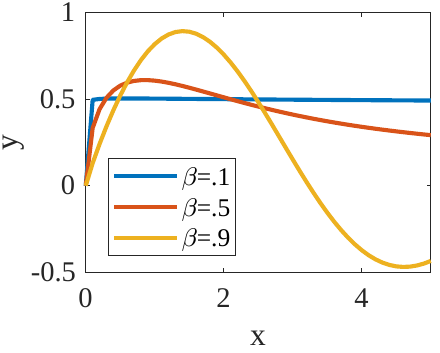}
\caption{A Monte Carlo simulation of the solution of $D_{c}^{\protect\beta %
}D_{c}^{\protect\beta }y(t)+\protect\omega ^{2}y(t)=\sin _{\protect\beta }{%
(t)}$ with $y(0)=0$ and $y^{\prime }(0)=0$.}
\label{fig4}
\end{figure}
\end{example}

\begin{example}
Let $0<\beta <1$, and consider 
\begin{equation*}
D_{c}^{\beta }D_{c}^{\beta }y(t)+\omega ^{2}y(t)=\frac{2t^{2\beta }}{\Gamma
(2\beta +1)}
\end{equation*}%
such that $y(0)=0$ and $y^{\prime }(0)=1$. It is simulated by solving $%
z^{\prime \prime }(t)+w^{2}z(t)=t^{2}$ such that $z(0)=0$ and $z^{\prime
}(0)=0$ using Runge-Kutta of 4 and 5 over simulated timelines $L(t)$ for each $t$.  Figure \ref{fig5} shows a Monte Carlo simulation of the solution; for more details see Section 5.

\begin{figure}[htb]
\centering
\includegraphics[scale=.4]{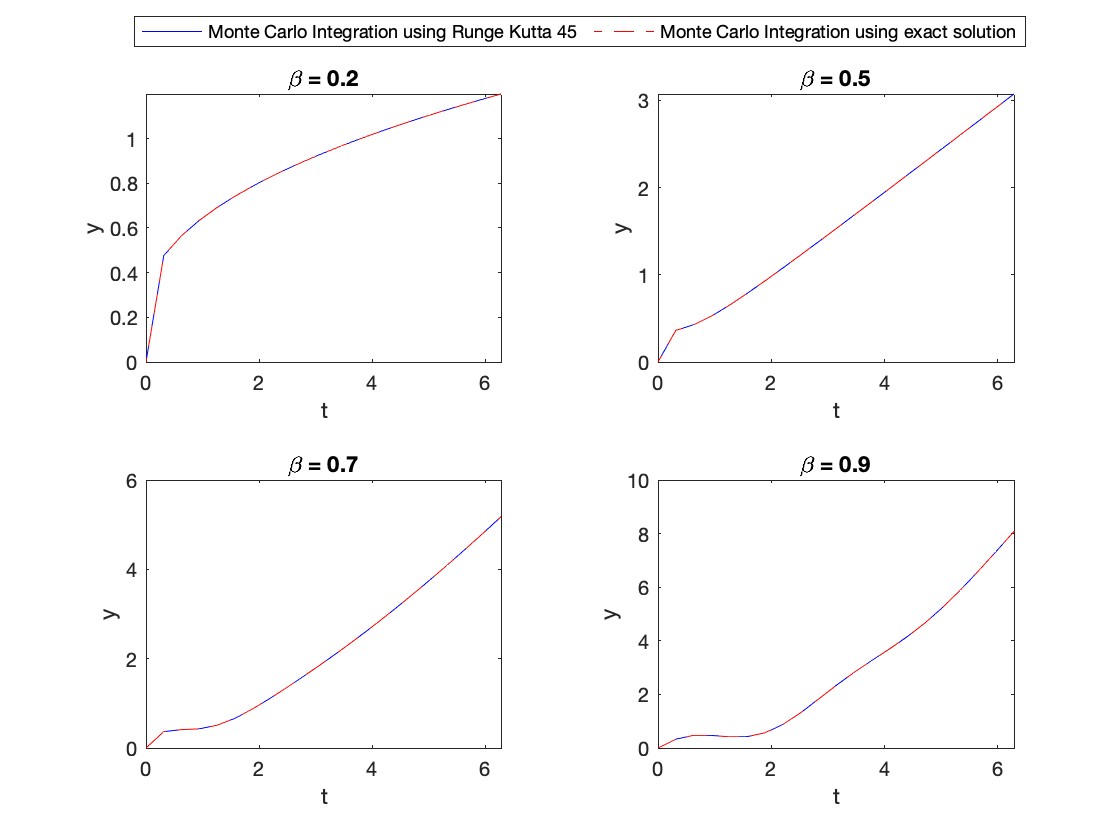}
\caption{A Monte Carlo simulation of the solution of $D_{c}^{\protect\beta %
}D_{c}^{\protect\beta }y(t)+\protect\omega ^{2}y(t)=\frac{2t^{2\protect\beta %
}}{\Gamma (2\protect\beta +1)}$ with $y(0)=0$ and $y^{\prime }(0)=1$ at $%
\protect\beta=0.2$ (top left), $\protect\beta=0.5$ (top right), $\protect%
\beta=0.7$ (bottom left), and $\protect\beta=0.9$ (bottom right).}
\label{fig5}
\end{figure}
\end{example}
\section{d'Alembert fractional formula for a fractional wave equation}

In this Section we present the last application of Theorem \ref%
{Theorem1} and Lemma \ref{lemma2}, the d'Alembert's formula for a fractional
wave equation. We also use the combination of Monte Carlo simulations with $g_{\beta }(\cdot ,t)$ to show how the solutions change on time  for different values of fractional derivatives, see Figure 9.

\begin{corollary}
(d'Alembert's formula)\ The fractional wave equation%
\begin{eqnarray}
D_{c}^{\beta }D_{c}^{\beta }u_{\beta }(x,t) &=&c^{2}\frac{\partial
^{2}u_{\beta }(x,t)}{\partial x^{2}},\text{ }-\infty <x<\infty ,\text{ }t>0 \label{Wave}
\\
u_{\beta }(x,0) &=&f_{\beta }(x) \\
\frac{\partial u_{\beta }(x,0)}{\partial t} &=&0
\end{eqnarray}%
has the following solution%
\begin{equation*}
u_{\beta }(x,t)=\frac{1}{2}[f_{\beta }(x+ct)+f_{\beta }(x-ct)].
\end{equation*}
\end{corollary}

Proof: Taking the Fourier transform on both sides of (6.1)-(6.3) and using the
representation of Caputo derivative presented in Lemma \ref{lemma2} , we
obtain%
\begin{eqnarray}
D_{c}^{\beta }\text{\ }D_{c}^{\beta }\mathcal{U}_{{}}(k,t) &=&-c^{2}k^{2}%
\mathcal{U}_{{}}(k,t)  \label{fractional fourier} \\
\mathcal{U}_{{}}(k,0) &=&F(k)  \notag \\
\text{\ }\frac{d^{{}}\mathcal{U}_{{}}}{dt^{{}}}(k,0) &=&0.  \notag
\end{eqnarray}%
(where $\mathcal{U}_{{}}(k,t)$ is the Fourier transform of $u(x,t)$) in order to solve
equation \ref{fractional fourier} using Theorem \ref{Theorem1}. We first
solve the regular ODE%
\begin{eqnarray*}
\text{\ }\frac{d^{2}\mathcal{U}}{dt^{2}}(k,t) &=&-c^{2}k^{2}\mathcal{U}(k,t)
\\
\mathcal{U}(k,0) &=&F(k) \\
\text{\ }\frac{d\mathcal{U}}{dt^{{}}}(k,0) &=&0
\end{eqnarray*}%
The solution is given by%
\begin{equation*}
\mathcal{U}(k,t)=\frac{F(k)}{2}\left( e^{ickt}+e^{-ickt}\right) .
\end{equation*}%
Taking the beta transformation, we obtain the solution given by%
\begin{eqnarray*}
\mathcal{U}_{\beta }(k,t) &=&\int_{0}^{\infty }g_{\beta }(z;t)%
\,\mathcal{U}(k,z)dz \\
&=&\int_{0}^{\infty }g_{\beta }(z;t)\frac{F(k)}{2}\left(
e^{ickz}+e^{-ickz}\right) dz.
\end{eqnarray*}%
Now, \ taking the Fourier inverse%
\begin{equation*}
\frac{1}{\sqrt{2\pi }}\int\limits_{%
\mathbb{R}
}e^{i\kappa x}\mathcal{U}_{\beta }(k,t)dk=\frac{1}{2}\frac{1%
}{\sqrt{2\pi }}\int\limits_{%
\mathbb{R}
}e^{i\kappa x}\int_{0}^{\infty }g_{\beta }(z;t)\frac{F(k)}{2}\left(
e^{ickz}+e^{-ickz}\right) dzdk
\end{equation*}%
and interchanging the order of integration%
\begin{equation*}
\frac{1}{\sqrt{2\pi }}\int_{0}^{\infty }g_{\beta }(z;t)\int\limits_{%
\mathbb{R}
}e^{i\kappa x}\mathcal{U}(k,z)dkdz=\frac{1}{2}\frac{1}{\sqrt{2\pi }}%
\int_{0}^{\infty }g_{\beta }(z;t)\int\limits_{%
\mathbb{R}
}e^{i\kappa x}\frac{F(k)}{2}\left( e^{ickz}+e^{-ickz}\right) dkdz,
\end{equation*}%
we get%
\begin{eqnarray*}
\int_{0}^{\infty }g_{\beta }(z;t)u(x,z)dz &=&\frac{1}{\sqrt{2\pi }}%
\int_{0}^{\infty }g_{\beta }(z;t)\int\limits_{%
\mathbb{R}
}e^{i(\kappa x+ckz)}\frac{F(k)}{2}dkdz \\
&&+\frac{1}{\sqrt{2\pi }}\int_{0}^{\infty }g_{\beta }(z;t)\int\limits_{%
\mathbb{R}
}e^{i(\kappa x-ckz)}\frac{F(k)}{2}dkdz.
\end{eqnarray*}%
Simplifying, we obtain
\begin{equation*}
u_{\beta }(x,t)=\frac{1}{2}\int_{0}^{\infty }g_{\beta }(z;t)f(x+cz)dz+\frac{1%
}{2}\int_{0}^{\infty }g_{\beta }(z;t)f(x+cz)dz
\end{equation*}%
Thus, we obtain the d'Alembert solution of the wave equation of the form%
\begin{equation*}
u_{\beta }(x,t)=\frac{1}{2}[f_{\beta }(x+ct)+f_{\beta }(x-ct)]=\frac{1}{2}%
\mathbb{E}\left[ f(x+\mathcal{T}_{\beta }(t))+f(x-\mathcal{T}_{\beta }(t))%
\right]
\end{equation*}
where the process $\mathcal{T}_{\beta }(t)$ can be seen as the inverse of a $%
\beta $ stable subordinator whose density is $g_{\beta }(x,t)$; for a
detailed discussion of the d'Alembert solution of the fractional wave
equation in general spaces, see \cite{Meerschaert2019} and references therein.
\begin{figure}[h]
{%
\includegraphics[scale=0.5]{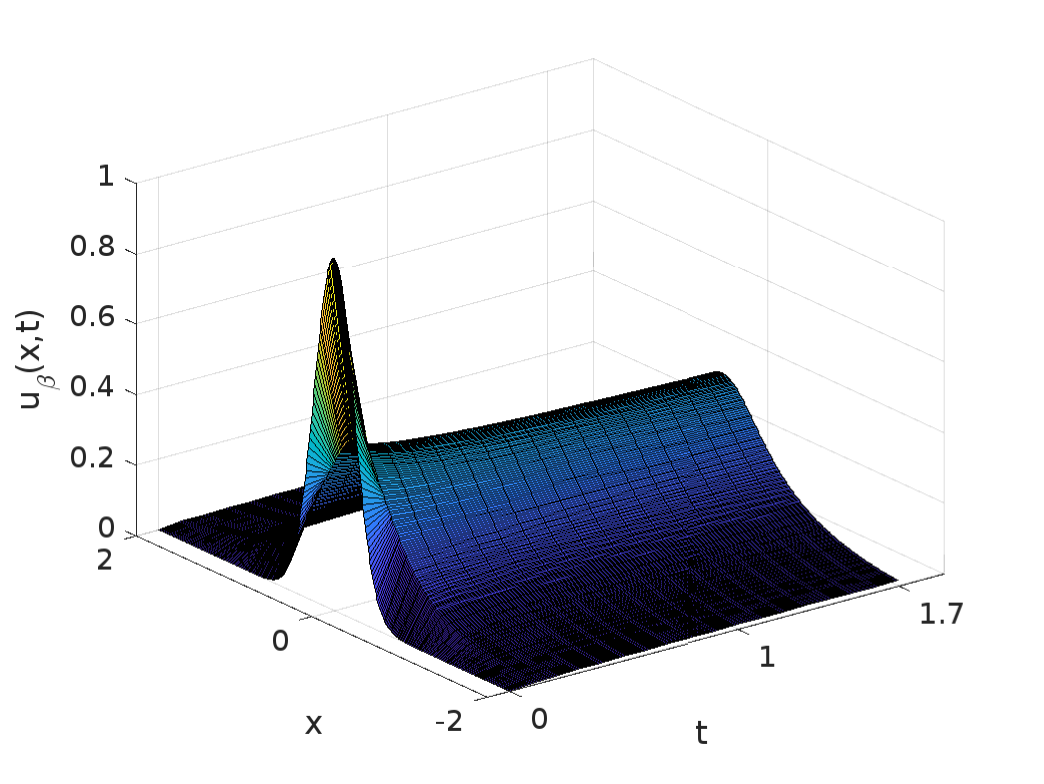}%
} {%
\includegraphics[scale=0.5]{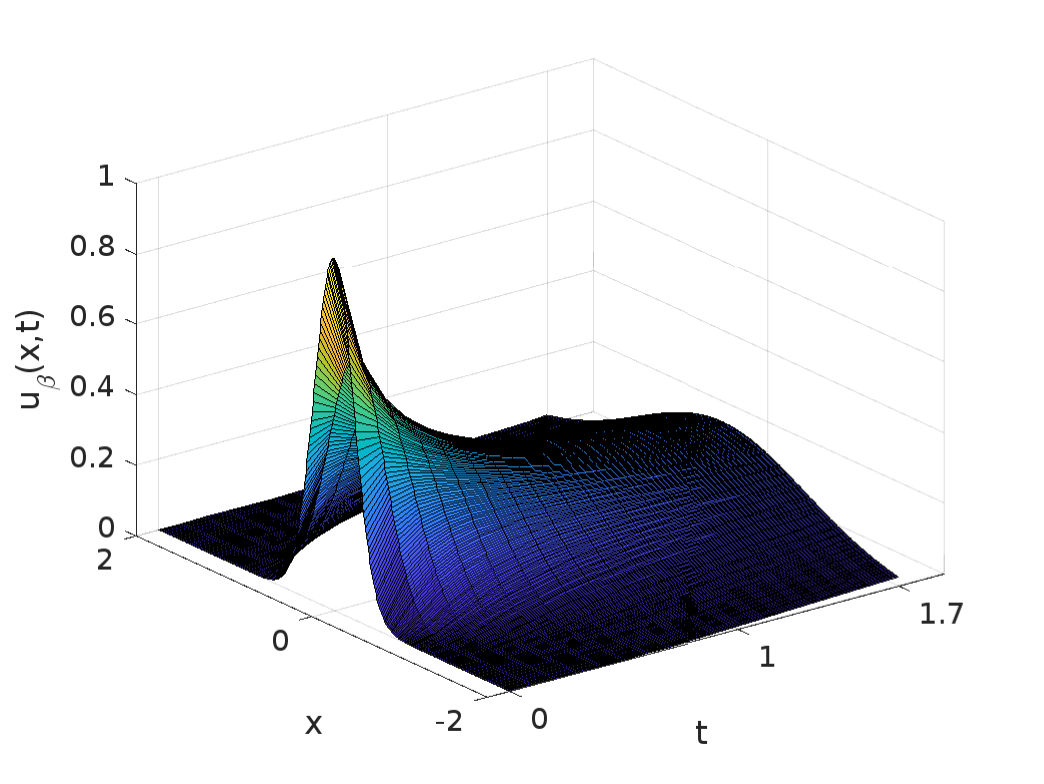}%
}

{%
\includegraphics[scale=0.5]{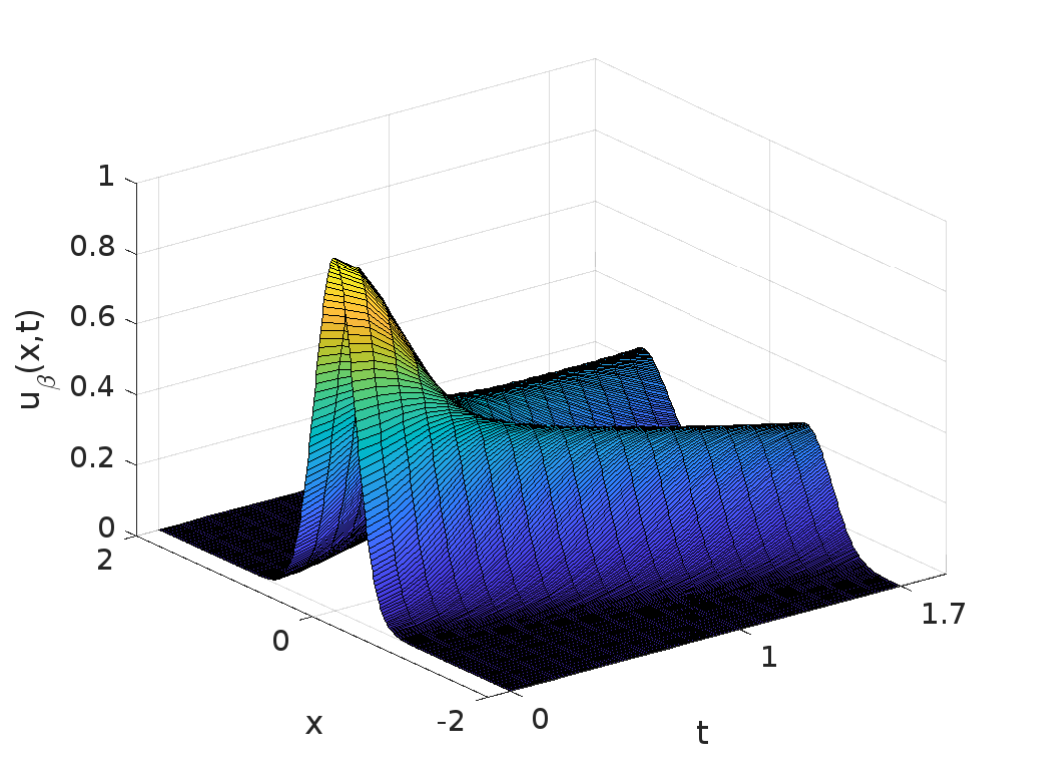}%
}
\caption{A Monte Carlo simulation of the solution for fractional wave equation $D_{c}^{\beta }D_{c}^{\beta }u_{\beta }(x,t) =c^{2}\frac{\partial
^{2}u_{\beta }(x,t)}{\partial x^{2}},
u_{\beta }(x,0) =f_{\beta }(x), 
\frac{\partial u_{\beta }(x,0)}{\partial t} =0$
 evaluated at $\protect%
\beta =0.1$ no diffusion (top), $\protect%
\beta =0.5$ showing diffusion (middle), and $\protect\beta =0.9$ 
 d'Alembert's formula (bottom) with $c=.5$ and
$f(x) = \exp(-10 x^2)$.}
\label{fig3}
\end{figure}

\section{Deep Learning to Simulate Fractional Differential Equations}

In this section, we will show how to use Feedforward Neural Networks
to simulate solutions of ordinary and fractional differential equations. The
main difference between the fractional and ordinary is in the number of past
measurements (memory length) required to predict the measurement of the next
step. We found that ODEs require one past step while FDEs require more than
one (longer memory).

The following algorithm is used to simulate a solution for the fractional
differential equations mentioned above. The algorithm uses Runge-Kutta with
hybrid orders 4 and 5 to solve the differential equations and then uses
the Monte Carlo interpolation to conduct the Wright-type transformation. The
data from the transform is used as training data for a feedforward neural network
(FFNN). The neural network then provides predictions on future solutions to
the fractional differential equation. Due to the definition of derivatives,
in fractional calculus being different than traditional derivatives we
approach the construction of the FFNN differently. For normal differential
equations we would be able to use a single moment from the past to predict a
moment in the future. For our fractional differential equation example we
needed to take three moments from the past to predict the future. The neural
network was constructed out of three layers with 10 neurons in each layer.

\textbf{Algorithm For Feedforward Neural Network}

\begin{enumerate}
\item \textbf{Step 1:} Simulate training, validation, and test Monte-Carlo
simulations of the fractional differential equations.

\item \textbf{Step 2:} Train the FFNN on the training data set and use the
validation data set.

\item \textbf{Step 3:} Use three data points or more to predict one future
point.
\end{enumerate}

\begin{example}
We show the FFNN solution of the FDE given by 
\begin{equation*}
D_{c}^{\beta }D_{c}^{\beta }D_{c}^{\beta }y(t)+2D_{c}^{\beta }D_{c}^{\beta
}y(t)+D_{c}^{\beta }y(t)+5y(t)=E _{\beta }{(-t)}
\end{equation*}
with $0<\beta<1$ with random $y(0)$.

\begin{figure}[htp]
{\includegraphics[scale=0.5]{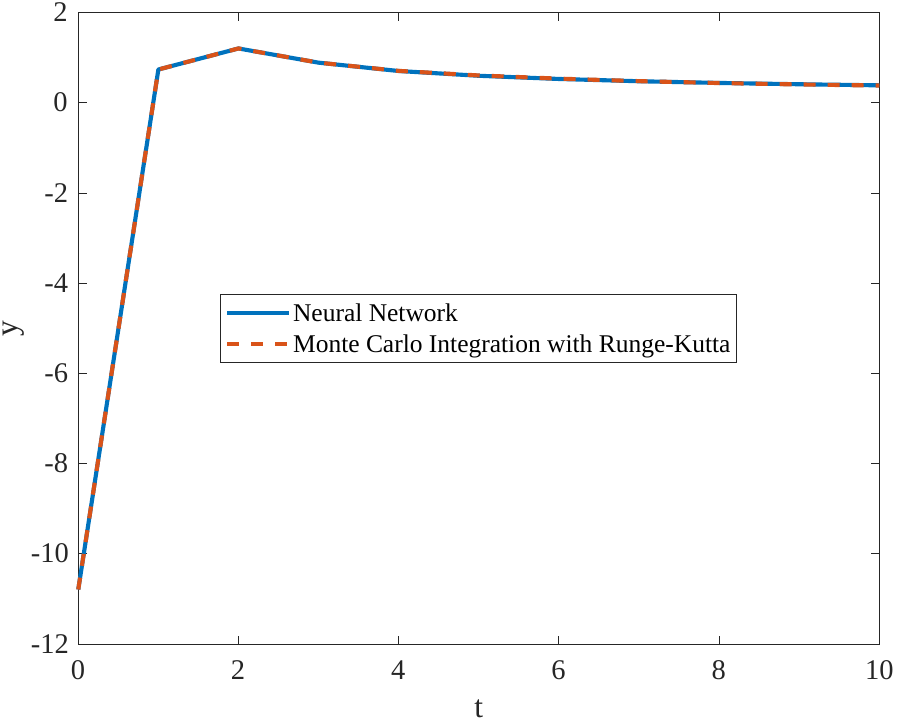}} {%
\includegraphics[scale=0.5]{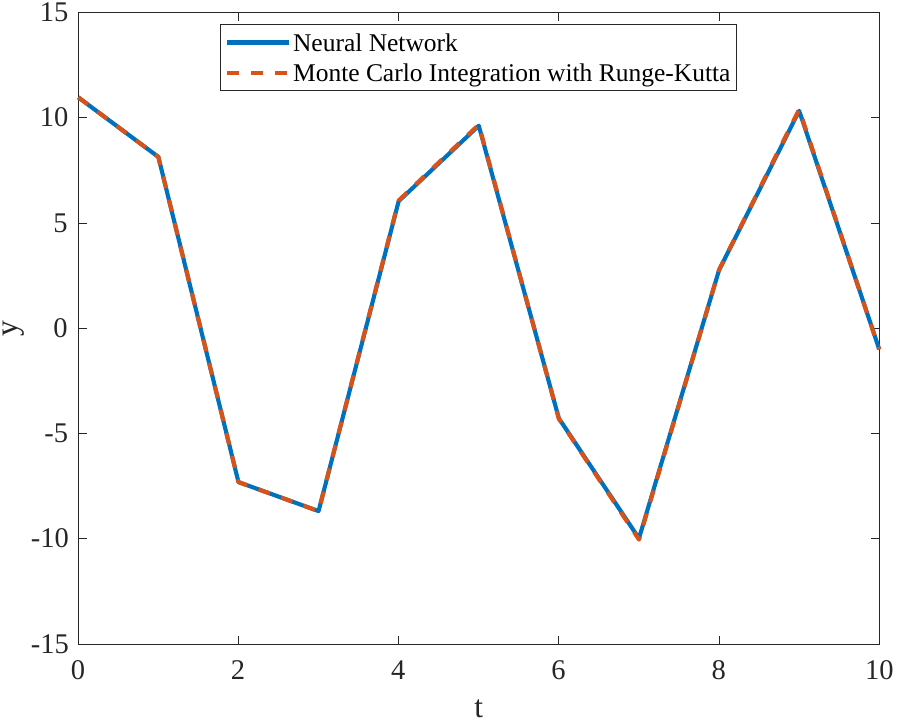}}
\caption{The graph shows neural network and Monte Carlo integration with
Runge-Kutta solutions for $D_{c}^{\protect\beta }D_{c}^{\protect\beta %
}D_{c}^{\protect\beta }y(t)+2D_{c}^{\protect\beta }D_{c}^{\protect\beta %
}y(t)+D_{c}^{\protect\beta }y(t)+5y(t)=E _{\protect\beta }{(-t^\protect\beta)%
}$ for $\protect\beta=0.5$ (left) and $\protect\beta=0.9$ (right).}
\label{fig1}
\end{figure}

The two graphs produced agree well. When $\beta=0.5,$ the mean squared error
between the neural network and the Monte Carlo simulations is $3.9833e-10$.
When $\beta=0.9$ the mean squared error between the neural network and the
Monte Carlo simulations is $2.0942e-9$. See Figure 9. 
\end{example}

\section{conclusion}
A new relationship (Lemma \ref{lemma2}) is established
between the solutions of higher fractional differential equations and a Wright-type transformation. Theorem \ref{Theorem1} allows us to solve fractional
differential equations of higher orders with certain initial conditions. Furthermore, solutions could be
interpreted as expected values of functions in a random-time process. Applications include the fractional beam equation, fractional electric
circuits with special functions as external sources such as the Mittag-Leffler
function and other special functions, and we provide a new derivation of d'Alembert's
formula for fractional wave equations.  
    
    We have shown that formulating solutions of
fractional differential equations and partial differential equations as
expected values with respect to power-law distributions could
be enabled using Monte Carlo integration methods along with other numerical methods such as Runge-Kutta to numerically solve fractional differential equations. The results of combining Monte Carlo integration methods along with Runge-Kutta methods showed an excellent approximation to the exact solution combined with Monte Carlo integration. Using a feedforward neural network along with a number of steps to make a prediction was also very successful in making similar simulation results like those done by Monte Carlo integration.

\begin{acknowledgement}
All the authors are honored and grateful to report that M. Nacianceno, Y.
Sepulveda, T. Stuck and J. Williams were fully  funded during their research experience for undergraduates (REU) by National
Science Foundation with Award \# 2150478 as undergraduate researchers. Also,
the faculty E. Suazo, T. Oraby, H. Rodrigo and J. Sifuentes were partially
funded by the same grant.
\end{acknowledgement}

\end{document}